\documentclass[10pt]{amsart}

\title{The information content of points on lines and $k$-plane extensions}

\author{Jacob B. Fiedler}
\address{Department of Mathematics, University of Wisconsin-Madison, Wisconsin 53715}

\email{jbfiedler2@wisc.edu}
\thanks{The author was supported in part by NSF DMS-2037851 and NSF DMS-2246906.}

\usepackage{amsmath}
\usepackage{amssymb}
\usepackage{mdframed}
\usepackage{amsthm}
\usepackage{mathtools}
\usepackage{enumitem}
\usepackage{placeins}
\usepackage{url}
\usepackage{hyperref}

\newtheorem{thm}{Theorem}

\newtheorem{lem}[thm]{Lemma}
\newtheorem{prop}[thm]{Proposition}

\theoremstyle{remark}

\DeclareMathOperator{\Dim}{Dim}

\newcommand{\R}{\mathbb{R}}

\newcommand{\N}{\mathbb{N}}
\newcommand{\Q}{\mathbb{Q}}

\newcommand{\ve}{\varepsilon}

\newenvironment{proofof}[1]{\begin{trivlist}
		\item[\hskip \labelsep \textit{Proof of #1.}]}{\end{trivlist}}

\begin{document}
\begin{abstract}
     We prove a new lower bound on the algorithmic information content of points lying on a line in $\mathbb{R}^n$. More precisely, we show that a typical point $z$ on any line $\ell$ satisfies 
    \begin{center}
        $K_r(z)\geq \frac{K_r(\ell)}{2} + r - o(r)$
    \end{center}
    at every precision $r$. In other words, a randomly chosen point on a line has (at least) half of the complexity of the line plus the complexity of its first coordinate. We apply this effective result to establish a classical bound on how much the Hausdorff dimension of a union of positive measure subsets of $k$-planes can increase when each subset is replaced with the entire $k$-plane. To prove the complexity bound, we modify a recent idea of Cholak-Cs\"ornyei-Lutz-Lutz-Mayordomo-Stull. 

\end{abstract}
\maketitle

\section{Introduction}

Two points determine a line, and approximations of those points determine an approximation of that line, provided that the points are suitably separated. Hence, it is reasonable to ask whether each point must contain at least half the amount of the information that the line contains. In general, this is not true; although a line in a random direction passing through the origin has a rather high complexity, it has a very simple point on it. However, we show that \emph{typical} points on a line actually satisfy a somewhat stronger inequality than the aforementioned half-complexity bound. In particular, we have the following as our main algorithmic result. 

\begin{thm}\label{thm:halfLines}
    Let $A\subseteq\mathbb{N}$, $\ell$ a line in $\mathbb{R}^n$, and $s\in \mathbb{R}$ be given. Assume that $s$ is random relative to $A$ and that 
    \begin{equation*}
    K^{A}_r(\ell \mid s)\geq K^A_r(\ell) - O(\log r). 
    \end{equation*}
    Then
    \begin{equation*}
        K^A_r(\ell(s))\geq \frac{K^A_r(\ell)}{2} + r - o(r).
    \end{equation*}
\end{thm}
\noindent Here and throughout the paper, $\ell(s)$ denotes the point on $\ell$ with $i$th coordinate $s$, where $i$ corresponds to the first coordinate axis onto which $\ell$ has a full projection, and $K^A_r$ indicates the Kolmogorov complexity at precision $r$, which is defined in Section 2. 

 The Kolmogorov complexity of points on lines is an interesting topic in algorithmic information theory, but via the point-to-set principle of J. Lutz and N. Lutz, it also has applications to very classical questions in geometric measure theory \cite{lutz2018algorithmic}. For instance, strong bounds on the complexity of points on lines in $\mathbb{R}^2$ imply estimates on the dimension of Furstenberg sets \cite{lutz2020bounding}. In this paper, we are able to use Theorem \ref{thm:halfLines} to prove the following result on $k$-plane extensions.

\begin{thm}\label{thm:classical}
Let $E\subseteq\mathbb{R}^n$ and let $F$ denote the union of $E$ with every $k$-plane that intersects $E$ in a set of positive measure. Then either $E=F$, or  

\begin{equation*}
    \dim_H(F)\leq 2 \dim_H(E) - k.
\end{equation*}
\end{thm}

Originally, this problem was posed for lineal extensions, that is, the case $k=1$. Keleti proved that the Hausdorff dimension of a collection of line segments in the plane does not increase when the line segments are replaced with the associated full lines and conjectured that this also holds for collections of line segments in $\mathbb{R}^n$ \cite{keleti2016lines}. In prior work with Bushling on line segment extensions, the author proved a strong packing dimension version of Keleti's theorem, and bounded how much the packing dimension can increase for sets of line segments in higher dimensions \cite{bushling2025extension}. The author recently generalized the latter bound to $k$-plane extensions \cite{fiedler2025extensionsUnions}, proving that
\begin{equation*}
    \dim_P(F)\leq 2 \dim_P(E) - k,
\end{equation*}
where $E$ and $F$ are as above. 

We briefly describe some other work relevant to the dimension of these extensions. Falconer and Mattila generalized Keleti's result to show that extensions of positive measure subsets of hyperplanes cannot increase the Hausdorff dimension of the set \cite{falconer2016strong}; this is stronger than Theorem 1 in the case $k=n-1$. Additionally, work of H\'era, Keleti, and M\'ath\'e implies that extending each full Hausdorff dimension subset of a hyperplane cannot increase the Hausdorff dimension of the set \cite{heraKeleti2019hausdorff}; the author showed that the packing dimension of the set also cannot increase under this extension \cite{fiedler2025extensionsUnions}. H\'era proved bounds on the Hausdorff dimension of Furstenberg-type sets, that is to say, unions of fractal subsets of $k$-planes \cite{Hera2019}. In the case $k=1$, H\'era's more general bound implies that a union of Hausdorff dimension $1$ subsets of lines satisfies
\begin{equation*}
    \dim_H(E)\geq 1 + \frac{t}{2},
\end{equation*}
where $t$ is the Hausdorff dimension of the set of lines (as a subset of the affine Grassmannian). If $F$ is the union of these lines, then we have the upper bound $\dim_H(F)\leq 1+t$. This implies 
\begin{equation*}
    \dim_H(F)\leq 2 \dim_H(E) - 1,
\end{equation*}
which is (a somewhat stronger version of) the $k=1$ case of Theorem \ref{thm:classical}. 

Finally, Keleti and M\'ath\'e showed that the Kakeya conjecture implies the line segment extension conjecture \cite{keleti2022equivalences}. Hence, Wang and Zahl's recent proof of the Kakeya conjecture in $\mathbb{R}^3$ implies that unions of line segments in $\mathbb{R}^3$ have the same Hausdorff dimension as unions of the corresponding full lines \cite{wang2025volumeestimatesunionsconvex}. It is worth noting that this implication also goes the other way, in part. Keleti showed that the line segment extension conjecture in $\mathbb{R}^n$ implies Besicovitch sets in $\mathbb{R}^n$ have Hausdorff dimension at least $n-1$, and that if this conjecture holds for \emph{all} $n$, then the packing dimension of Besicovitch sets is $n$ \cite{keleti2016lines}. 

\subsection{Overview of the proofs}

We briefly describe some of the ideas of the proofs. Theorem \ref{thm:classical} follows from Theorem \ref{thm:halfLines} using the point-to-set principle and a radial slicing argument. We want to upper bound the information content of a point $x\in F$ lying on some line $\ell$ to which we apply Theorem \ref{thm:halfLines}, but we have enough control over that line that it can be forced to contain a significant amount of information unrelated to $x$. Since $x$ and this extra information can be determined from $\ell$ and one coordinate of $x$, and $\ell$ in turn cannot be too much more complex than typical points in $\ell \cap E$, we have the desired bound. 

The proof of Theorem \ref{thm:halfLines} is a bit harder, and utilizes a very interesting new technique of Cholak-Cs\"ornyei-Lutz-Lutz-Mayordomo-Stull \cite{CholakCsorn2025Bourgain}. The problem to overcome is that although a simple argument gives,
\begin{equation*}
    K_r^A(\ell)\leq K^A_r(y, z) + O(\log r) \leq K^A_r(y) + K^A_r(z) + O(\log r),
\end{equation*}
where $y$ and $z$ are appropriately separated points on $\ell$, it is not obvious that $y$ and $z$ have comparable complexity at every precision. One could imagine that $y$ has high complexity whenever $z$ has low complexity, and vice versa, which is consistent with $\ell$ having high complexity but prevents either $y$ or $z$ from satisfying the conclusion of Theorem \ref{thm:halfLines}.  

The idea of \cite{CholakCsorn2025Bourgain} is that to overcome an obstruction like this, one can work with surrogate points. These authors apply their technique to generalize Bourgain's exceptional set estimate for projections \cite{Bourgain10} and to improve bounds on the dimension of pinned distance sets in the plane. In the context of points on a line, if the conclusion of Theorem \ref{thm:halfLines} fails for some $\ell$ and $s$, then once can show that there are ``many'' values $s^\prime$ and lines $\ell^\prime$ such that $\ell^\prime(s^\prime)$ has lower than the expected complexity. A combinatorial argument then allows one to find a pair $u, v\in \mathbb{R}$ that supports many such lines -- in particular a large enough number that one of these lines $l$ has complexity higher than that of the original line $\ell$. On one hand, $l(u)$ and $l(v)$ determine $l$ which has high complexity compared to $\ell$. On the other hand, $l(u)$ and $l(v)$ both have, by assumption, low complexity compared to $\ell$ at precision $r$, which yields a contradiction. 

As compared to \cite{CholakCsorn2025Bourgain}, in this paper we actually ask somewhat more of the pair $(u, v)$, namely that it also determines a significant amount of information unrelated to $\ell$, which is responsible for Theorem \ref{thm:halfLines}'s extra ``$+r$''. This necessitates modifications to the proof, but it is possible because we are able to make a strong assumption on the complexity of $s$. 

As a final comment, we note that the surrogate point technique of \cite{CholakCsorn2025Bourgain} is an important extension of the algorithmic toolkit as it applies to geometric measure theory. Previously, many new dimension bounds in GMT arising from effective methods required the application of one or more ``enumeration'' lemmas. These lemmas tend to work well when one has a great deal of control over some of the objects in the problem. For instance, some of the strongest bounds in \cite{fiedler2025extensionsUnions} are for unions of large subsets of hyperplanes, as one can choose a point in such a subset to have complexity about $(n-1) r$ relative to any given object; for other examples see \cite{lutz2020bounding, stull2022pinneddistancesetsusing}. However, this surrogate point technique applies in some situations where very \emph{little} complexity exists to take advantage of, such as projections onto directions with low effective Hausdorff dimension, which can be compared to prior results on projections that employed enumeration lemmas \cite{LutStu18Projections, fiedler2025universalsetsprojections}. 

\section{Preliminaries}

Let $K^A(\sigma)$ denote the (prefix-free) Kolmogorov complexity of a finite-data object $\sigma$ (relative to an oracle $A\subseteq\mathbb{N}$), and let $K^A(\sigma\mid\tau)$ denote the Kolmogorov complexity of $\sigma$ given $\tau$\footnote{We recommend \cite{downey2010} to readers unfamiliar with Kolmogorov complexity and \cite{lutz2018algorithmic} for more details on the definitions we briefly introduce below.}. Define the Kolmogorov complexity at precision $r$ of $x\in\mathbb{R}^n$ as
\begin{equation*}
K^A_r(x)=\min\left\{K(p)\,:\,p\in B_{2^{-r}}(x)\cap\Q^m\right\},
\end{equation*}
and the Kolmogorov complexity of $x$ at precision $r$ given $y$ at precision $s$ as 
\begin{equation*}
K^A_{r,s}(x\mid y)=\max\big\{\min\{K_r(p\mid q)\,:\,p\in B_{2^{-r}}(x)\cap\Q^m\}\,:\,q\in B_{2^{-s}}(y)\cap\Q^n\big\}.
\end{equation*}
In prior work \cite{fiedler2025extensionsUnions}, using ``rational'' projection matrices, the author similarly defined the complexity at precision $r$ of elements of the Grassmannian and affine Grassmannian. 

In Euclidean space, we can associate a $2^{-r}$ dyadic rational to any point $x$ by, say, rounding down in each coordinate. The structure of the Grassmannian makes it somewhat less obvious what the dyadics ``should'' be, but \cite{fiedler2025extensionsUnions} defined subsets of the rational elements with enough relevant properties to assume the same role as the dyadics in these computability-theoretic arguments. Hence, let $d(x, r)$ denote the $2^{-r}$-``dyadic'' associated to $x$ (where $x$ is an element of the space $\mathcal{X}$, which is either $\mathbb{R}^n$, $\mathcal{G}(n, k)$, or $\mathcal{A}(n, k)$). Furthermore let $\mathcal{D}(\mathcal{X}, r)$ denote the set of dyadic elements of the space $\mathcal{X}$. An important fact is that the ``minimizing rational'' definition of the Kolmogorov complexity at precision $r$ of a point agrees with the complexity of the associated $2^{-r}$-dyadic up to a small error term \cite{lutz2020bounding, fiedler2025extensionsUnions}. In particular, there exists some $C_{\mathcal{X}}$ such that for every $x\in\mathcal{X}$
\begin{equation}\label{eq:approximationsMain}
    K^A(d(x, r))- C_{\mathcal{X}} \log r \leq K^A_r(x)\leq K^A(d(x, r))+ C_{\mathcal{X}} \log r.
\end{equation}
A similar statement holds for conditional complexity. 

We will repeatedly make use of a basic fact regarding the Kolmogorov complexity of elements of a computably enumerable set of dyadics (see, for instance, \cite{csornyei2025improvedboundsradialprojections}). We let $\sigma \in D_r\subseteq \mathcal{D}(\mathcal{X}, r)$, and we work relative to an oracle $A$. If there is a fixed Turing machine $M$ that enumerates the elements of $D_r$ given an input string $\tau$, then we can compute $\sigma$ given $\tau$ and $\sigma$'s place in $M$'s enumeration of $D_r$. The latter is an integer no larger than $\vert D_r\vert $, so
\begin{equation*}
    K^A(\sigma) \leq \log\vert D_r\vert + \vert \tau \vert + O(\log \log \vert D_r\vert + \log \vert \tau\vert).
\end{equation*}
\noindent In particular, if $M$ generates $D_r$ from $r$ and at most a logarithmic amount of extra information, and $D_r$ always contains the $2^{-r}$-dyadic associated to $x$, then
\begin{equation*}
    K_r^A(x) \leq \log\vert D_r\vert + O(\log r).
\end{equation*}

Kolmogorov complexity also obeys the symmetry of information. N. Lutz and Stull proved this in Euclidean space \cite{lutz2020bounding}, and recent work extended it to the Grassmannian and affine Grassmannian \cite{fiedler2025extensionsUnions}. 
\begin{prop}\label{prop:symmetry}
    Let $x\in\mathcal{X}$ and $y\in\mathcal{Y}$, where these spaces are as above. For every $A\subseteq\mathbb{N}$ and precisions $r, s$, 
    \begin{equation*}
        K^A_{r, t}(x, y) = K^A_t(y) + K_{r, t}^A(x\mid y) \pm O(\log (r + t)).
    \end{equation*}
\end{prop}

We mention a few more facts that will be useful in our algorithmic arguments. In general, oracle access to extra information is essentially no less useful than precision $s$ access, which is expressed in the following inequality:
\begin{equation*}
    K_r^{A, y}(x)\leq K_{r, s}^A(x\mid y) + O(\log(r + s)) \leq K_r^A(x) + O(\log (r + s)).
\end{equation*}

Additionally, one can bound how much the Kolmogorov complexity of an object changes when we alter the precisions. In \cite{case2015dimension}, Case and J. Lutz proved that,
\begin{prop}\label{prop:caseLutz}
For any $A \subseteq \mathbb{N}$, $r,t \in \N$, and $x \in \R^n$,
\begin{equation*}
    K_{r+t}^A(x) \leq K_r^A(x) + nt + O(\log(r + t)).
\end{equation*}
\end{prop}
Note that this Proposition is stated in Euclidean space; although it is very likely that such a bound also holds for elements of the Grassmannian, to our knowledge, a full proof of this fact has not been written. However, Proposition 25 in \cite{fiedler2025extensionsUnions} established that 
\begin{equation*}
    K_r^A(\ell) = K_r^A(a_1, ... a_{n-1}, b_1, b_{n-1}) \pm O(\log r),
\end{equation*}
where $a_i$ and $b_i$ are such that $\ell(s) = (s, a_1s + b_1, ..., a_{n-1} s+ b_{n-1})$. 

Since $(a_1, ... a_{n-1}, b_1, b_{n-1})\in \mathbb{R}^{2(n-1)}$ and the complexity of points and lines is unaffected by computable rotations (up to a logarithmic term), we essentially have the Case-Lutz bound for every line $\ell$, that is
\begin{equation}\label{eq:pseudoCaseLutz}
    K^A_{r+t}(\ell) \leq K^A_r(\ell) + 2(n-1) t + O(\log (r+t)). 
\end{equation}

There are notions of effective fractal dimension due to J. Lutz \cite{Lutz03a, Lutz03b}, and Mayordomo proved that these notions are equivalent to the following quantities \cite{Mayordomo02}, which we take to be the \emph{definition} of effective dimension,
\begin{equation*}
    \dim^A(x):= \liminf_{r\to\infty} \frac{K_r^A(x)}{r} \qquad \text{and} \qquad \Dim^A(x):= \limsup_{r\to\infty}  \frac{K_r^A(x)}{r}.
\end{equation*}

Finally, we can state the point-to-set principle of J. Lutz and N. Lutz \cite{lutz2018algorithmic}. For any $E\subseteq\mathbb{R}^n$, 
\begin{equation*}
    \dim_H(E)=\min_{A\subseteq\mathbb{N}}\sup_{x\in E}\dim^A(x) \qquad \text{and} \qquad \dim_P(E)=\min_{A\subseteq\mathbb{N}}\sup_{x\in E}\Dim^A(x).
\end{equation*}
Point-to-set principles also hold for subsets of the Grassmannian and the affine Grassmannian \cite{LuLuMay2023PtS, fiedler2025extensionsUnions}.

\section{Bounding the Hausdorff dimension of $k$-plane extensions}

We will begin by deriving Theorem \ref{thm:classical} (restated below) from Theorem \ref{thm:halfLines}. 

\vspace{2mm}

\noindent \textbf{Theorem \ref{thm:classical}.} \textit{Let $E\subseteq\mathbb{R}^n$ and let $F$ denote the union of $E$ with every $k$-plane that intersects $E$ in a set of positive measure. Then either $E=F$, or } 

\begin{equation*}
    \dim_H(F)\leq 2 \dim_H(E) - k.
\end{equation*}

\begin{proof}
    Using the point-to-set principle, it suffices to show\footnote{c.f. Theorem 3's proof in \cite{fiedler2025extensionsUnions}.} that for every oracle $A$, $V\in\mathcal{A}(n, k)$, and positive measure $I\subseteq V$,
    \begin{equation*}
 \sup_{x\in V} \dim^A(x)\leq 2 \sup_{z\in I}\dim^A(z) - k.
    \end{equation*}
    Let $A, V, I$ and $x\in V$ be given. For every $x$ and $\ve>0$, we will find some $z$ such that $\dim^A(x)\leq 2 \dim^A(z) - k + \ve$. Note that we are finding a single such $z\in I$ instead of a pair of points $y, z$. Hence, it suffices to show that $z$ satisfies 
    \begin{equation}\label{eq:reducedBound}
        K^A_r(x)\leq 2 K^A_r(z) - kr + \ve r 
    \end{equation}
   for all $r$ sufficiently large. 
   
   We assumed that $I$ has positive $k$-dimensional measure, so by Fubini's theorem, there is a positive $(k-1)$-dimensional measure set of lines through $x$ each intersecting $I$ in a set of positive $1$-dimensional measure. Choose $\ell$ from this set such that
  \begin{equation}\label{eq:highComplexityLine}
  K_r^{A, x}(\ell)\geq (k - 1) r - O(\log r).
  \end{equation}
  Morally speaking, this is responsible for the extra ``$-k$'' term in the theorem. 
  
  Because computable rotations only affect complexity by at most a logarithmic term, we may assume without loss of generality that the projection of $\ell\cap I$ onto the first coordinate has positive measure. Hence, we may choose $s$ in this projection to be random relative to $A$ and $\ell$. In particular, it is random relative to $A$, so the first condition of Theorem \ref{thm:halfLines} holds. To see that the second holds, note that by our choice of $s$ and the symmetry of information,
  \begin{align*}
      K^A_r(s) + K^A_r(\ell) &\leq K^{A, \ell}_r(s)  + K^A_r(\ell)+ O(\log r)\\ 
      &\leq K^A_r(s\mid \ell) + K^A_r(\ell)+ O(\log r) \\
      &= K_r^A(s, \ell)+ O(\log r)\\
      &= K_r^A(\ell\mid s) + K_r^A(s)+ O(\log r),
  \end{align*}
  so
  \begin{equation*}
      K^A_r(\ell\mid s)\geq  K^A_r(\ell) - O(\log r).
  \end{equation*}

   We define $z= \ell(s)$. Note that a precision $r$ approximation of $x$ is computable from a precision $r$ approximation of $\ell$ and at most $r + O(\log r)$ bits of additional information (see, for instance, Lemma 24 of \cite{fiedler2025extensionsUnions}). Hence, using the symmetry of information,
   \begin{align*}
       K^A_r(\ell) + r &\geq K^A_r(x, \ell) - O(\log r)\\
       &= K^A_r(\ell \mid x) + K^A_r(x) - O(\log r)\\
        &\geq K^{A, x}_r(\ell) + K^A_r(x) - O(\log r)\\ 
        &\geq  (k-1) r + K^A_r(x) - O(\log r).
   \end{align*}

By our choice of $z$, we can apply Theorem \ref{thm:halfLines}. For sufficiently large $r$, in combination with the above, this yields
\begin{align*}
    K^A_r(x) &\leq K^A_r(\ell) - (k-2) r + O(\log r)\\
        &\leq 2K_r^A(z) - 2r - (k-2) r + \ve r + O(\log r)\\
        &\leq 2 K^A_r(z) - k r + \ve r
\end{align*}
which establishes \eqref{eq:reducedBound} and completes the proof. 

\end{proof}

\section{Proof of Theorem \ref{thm:halfLines}}

Now, we turn to the proof of Theorem \ref{thm:halfLines}. This proof will require a key combinatorial lemma,  which is Lemma 1 in \cite{CholakCsorn2025Bourgain}.

\begin{lem}\label{lem:combinatorial}
       Let $L$ and $V$ be finite sets, $\alpha\in(0,1)$, and, for each $d\in L$, $\sim_d$ a binary relation on $V$. Suppose there is a collection $\{N_v\mid v\in V\}$ of finite subsets of $L$ such that, for all $v\in V$, we have $|N_v|\geq\alpha\vert L\vert$ and, for all $d\in N_v$,
    \begin{equation*}\label{eq:fewsim}
        |\{u\in V\mid u\sim_d v\}|\leq \frac{\alpha^2}{4}|V|.
    \end{equation*}
    If there is a single binary relation $\sim$ such that $\sim_d$ is $\sim$ for all $d\in L$, then there exist $u,v\in V$ such that $u\not\sim v$ and $|N_u\cap N_v|>\frac{\alpha^2}{2}\vert L\vert $.
\end{lem}

We restate Theorem \ref{thm:halfLines} for convenience.

\noindent \textbf{Theorem \ref{thm:halfLines}.} Let $A\subseteq\mathbb{N}$, $\ell\in \mathcal{A}(n, 1)$, and $s\in \mathbb{R}$. Assume that $s$ is random relative to $A$ and that 
    \begin{equation*}
    K^{A}_r(\ell \mid s)\geq K^A_r(\ell) - O(\log r). 
    \end{equation*}
    Then
    \begin{equation*}
        K^A_r(\ell(s))\geq \frac{K^A_r(\ell)}{2} + r - o(r).
    \end{equation*}

Before we begin the proof, we make a few comments. It is reasonable to compare this theorem to Theorem 3 in \cite{CholakCsorn2025Bourgain}, which bounds the complexity of orthogonal projections. The geometry is similar; a line in $\mathbb{R}^n$ can be approximated using two separated points on it, and a point in $\mathbb{R}^2$ can be approximated using its projections onto two known, separated directions. We opt to prove Theorem \ref{thm:halfLines} at precision $r$ instead of on intervals of precisions $[r_1, r_2]$ as in \cite{CholakCsorn2025Bourgain}, which streamlines the latter portion of the proof. Additionally, in our application of Theorem \ref{thm:halfLines}, we do not need a bound on the quantity $K^{A, s}(\ell (s))$, so we can avoid a few steps pertinent to the extra oracle. On the other hand, our lower bound is larger by $r$, which we achieve by choosing the surrogate pair $(u, v)$ to satisfy an extra complexity bound. This requires a refinement of one of the relevant sets and additional work with the binary relation $\sim$ to which we apply Lemma \ref{lem:combinatorial}.

Finally, we have to name a number of constants $C_i$ in the proof, as we need to define several sets and Turing machines precisely. This introduces some clutter, but readers are advised that most of these are just instances of standard $O(\log r)$ error terms. 

\begin{proofof}{Theorem \ref{thm:halfLines}}

We may assume without loss of generality that the distance between $\ell$ and the line given by the first coordinate axis (call this axis $\ell_x$) is less than $\frac{1}{2}$, since rational rotations and translations only affect the complexity of points by at most a logarithmic term. Let $\ve>0$ and $r$ sufficiently large be given and assume for the purpose of contradiction that $K^A_r(\ell(s))< \frac{K^A_r(\ell)}{2} + r - \ve r$. 

\noindent \textbf{Some definitions and observations:} Set $t=\frac{\ve r}{2 n}$. Define 

\begin{equation*}
    L = \{d\in \mathcal{D}(\mathcal{A}(n, 1), r)\cap B_{1}(\ell_x): K^A(d)\leq K_r(\ell)+C_1 \log r\}
\end{equation*}
and
\begin{equation*}
    S = \{d\in \mathcal{D}(\mathbb{R}, r)\cap [0, 1]: K^A(d)\leq K^A_r(s) + C_2 \log r\text{ and } K^A_t(d)\leq K^A_t(s)+ 2 C_2 \log t\},
\end{equation*}
where $C_1$ and $C_2$ are the constants (for the appropriate space) in \eqref{eq:approximationsMain}. The extra conditions on \emph{where} these dyadics are located guarantee that the sets are finite. Furthermore, they ensure that if $\ell_1, \ell_2\in L$ and $s_1, s_2\in S$ are such that 
\begin{equation*}
\vert \ell_1 - \ell_2 \vert<\delta \qquad \text{ and } \qquad \vert s_1 - s_2 \vert<\delta, 
\end{equation*}
then
\begin{equation}\label{eq:smallSlopes}
    \vert \ell_1(s_1) - \ell_2(s_2) \vert \leq c\delta 
\end{equation}
for some universal constant $c$.

It is straightforward to adapt the proof of the first claim in Appendix C of \cite{CholakCsorn2025Bourgain} using Section 3.3 of \cite{fiedler2025extensionsUnions}, establishing that
\begin{equation}\label{eq:LIsLarge}
 2^{K^A_r(\ell) - O(\log r)} \leq \vert L\vert \leq  2^{K_r^A(\ell)+2C_1 \log r + 1}.
\end{equation}

Next, we define a collection of subsets of $L$ depending on a first coordinate $v$,

\begin{equation*}
    N_v = \{d\in L: K^A_r(d(v))< \frac{K^A_r(\ell)}{2} + r - \ve r + C_5 \log r\}.
\end{equation*}

The intuition is that for a given $v$, we want to describe the lines $d$ are such that the $d(v)$ has lower than the expected complexity. We will be able to use the assumption at the start of this proof (which essentially says that $\ell \in N_s$) to show that ``many'' $N_v$ have a large cardinality; this is essentially the hypothesis of Lemma \ref{lem:combinatorial}. To do so, we will need to make a few more preliminary definitions that will guarantee the set $V$ of these desirable $v$ has the appropriate properties. First, let
\begin{equation*}
    W_1 = \{d \in S: \vert N_{d}\vert > 2^{K^A_r(\ell) - C_3 \log r}\}.
\end{equation*}
Given $r, t, K^A_r(\ell), K^A_t(s)$, and $K^A_r(s)$, the sets $L$ and $S$ are computably enumerable.  Hence, assuming a Turing machine also has access to $\ve$, $W_1$ is computably enumerable. Note that we can choose $r$ sufficiently large depending on $\ve$, and we have that $t<r$, so all of the above are at most $O(\log r)$. Using the computable enumerability of $W_1$ and the above fact, we can define a Turing machine $M$ that, given $O(\log r)$ bits of information, 
\begin{enumerate}
    \item Lists all the elements $d\in \mathcal{D}(\mathbb{R}, t)$,
    \item Internally enumerates the $d^\prime\in W_1$
    \item Associates each $d^\prime$ to the $d$ such that $\lfloor d^\prime \rfloor_t = d$
    \item Outputs all the $d^\prime$ associated to a given $d$ as soon as there are $2^{K_r^A(s) - K_t^A(s) - C_4 \log r}$ such $d^\prime$; such a $d$ is called ``crowded''. 
\end{enumerate}
  Define $W$ to be the set of $d^\prime$ that $M$ outputs, and define $W^\prime$ to be the set of crowded $2^{-t}$-dyadics. Note that $W\subseteq W_1$. Our next goal is to bound the size of $W$. In particular, we need to show
  \begin{equation}\label{eq:boundOnW} 
  \vert W \vert \geq 2^{K^A_r(s) - O(\log r)}
  \end{equation}
  Note that 
  \begin{equation*}
      \vert W \vert = 2^{K^A_r(s) - K^A_t(s)- C_4 \log r} \vert W^\prime \vert, 
  \end{equation*}
   by definition. Hence, \eqref{eq:boundOnW} follows if
   \begin{equation}\label{eq:boundOnWPrime}
       \vert W^\prime \vert\geq 2^{K^A_t(s) - O(\log r)}.
   \end{equation}
 
  \noindent \textbf{Bounding the size of $W^\prime$:} Briefly, we know $W^\prime$ is large because $d(s, t)$\footnote{Recall the notation for dyadic elements.} is in $W^\prime$ and $W^\prime$ is computably enumerable, meaning we can bound its size in terms of the complexity of $s$ at precision $t$. However, showing that $d(s, t)\in W^\prime$ requires several steps. 
  
  For ease of reference, let $s_r:=d(s, r)$. First, we show that $d(\ell, r)$ is in $N_{s_r}$. Clearly, the definition of $C_1$ ensures $d(\ell, r)\in L$. An application of \eqref{eq:smallSlopes} allows us to see $\vert d(\ell, r)(s_r) - \ell(s_r)\vert < c 2^{-r}$. In conjunction with Proposition \ref{prop:caseLutz} and \eqref{eq:approximationsMain}, the above bound implies that there is some constant $C_5$ such that 
  \begin{equation*}
  \vert K^A(d(\ell, r)(s_r)) - K^A_r(\ell(s))\vert \leq C_5 \log r.
  \end{equation*}
  Combining our assumption on $K_r^A(\ell(s))$ at the start of the proof with the above gives
  \begin{equation*}
      K^A(d(\ell, r)(s_r)) \leq \frac{K^A_r(\ell)}{2} + r - \ve r + C_5 \log r,
  \end{equation*}
  which was exactly the condition for the inclusion of an element of $L$ in $N_{s_r}$.

In the following, we employ the fact that $d(\ell, r)\in N_{s_r}$, the fact that $N_{s_r}$ is computably enumerable given $s_r$ and little additional information, and the condition in Theorem \ref{thm:halfLines} that precision $r$ access to $s$ is unhelpful in computing $\ell$. In particular,
  \begin{align*}
      \log (\vert N_{s_r} \vert)
      &\geq K^{A}(d(\ell, r)\mid s_r) - O(\log r)\\ 
    &= K^{A}_r(\ell\mid s) - O(\log r)\\
    &\geq K^A_r(\ell) - O(\log r).
  \end{align*}
Hence, there exists some constant $C_3$ such that
\begin{equation}\label{eq:sInVOne}
    \vert N_{s_r}\vert > 2^{K_r^A(\ell) - C_3 \log r}.
\end{equation}
\noindent Recalling the definition of $W_1$, we see that $s_r\in W_1$. 

Next, we would like to use this fact to show that many nearby dyadics are \emph{also} in $W_1$, which will imply that $d(s, t)\in W^\prime$, i.e. that $d(s, t)$ is crowded. Let $T$ denote the set of $w\in W_1$ such that $d(w, t)= d(s, t)$. 

Again, we leverage computable enumerability; since $W_1$ is c.e. given a logarithmic amount of information, we have that $T$ is c.e. given $d(s, t)$  and a logarithmic amount of information. Hence, 
\begin{align*}
    \log(\vert T \vert) &\geq K^A(s_r\mid d(s, t)) - O(\log r)\\
    &= K_{r, t}^A(s) - O(\log r)\\
    &= K_{r}^A(s) - K_t^A(s)- O(\log r).
\end{align*}
So for some constant $C_4$, 
\begin{equation}\label{eq:sIsCrowded}
    \vert T \vert > 2^{K_r^A(s) - K^A_r(s) - C_4 \log r},
\end{equation}
which by definition means $d(s, t)$ will be labeled crowded by $M$. Once again this will allow us to exploit that a set -- in this instance $W^\prime$ -- is computably enumerable.

Since $W^\prime$ is c.e. given at most a logarithmic amount of information, and $d(s, t)\in W^\prime$, 
\begin{align*}
    \vert W^\prime\vert &\geq K^A(d(s, t)) - O(\log r)\\
     &= K^A_t(s)  - O(\log r),
\end{align*}
which establishes \eqref{eq:boundOnWPrime} and hence \eqref{eq:boundOnW}.

\noindent \textbf{Applying Lemma \ref{lem:combinatorial}:} We are now in a position where we can apply Lemma \ref{lem:combinatorial}. However, we will need to slightly refine the set $W$ so that the elements that come out of this application all have a high complexity. Hence, we define
\begin{equation*}
V =  \{v \in W: K^A_r(v)\geq r - C_6(\log r)\},
\end{equation*}
for $C_6$ large enough to guarantee that $\vert V\vert\geq \frac{\vert W \vert}{2}$. 

It is straightforward to see that there exists such a $C_6$ independent of $r$. Our assumption that $s$ was random relative to $A$ combined with \eqref{eq:boundOnW} ensures that
\begin{equation*}
    \vert W\vert\geq  2^{r - O_s(\log r)}.
\end{equation*}
The conclusion follows from the definition of Kolmogorov complexity; because there are fewer than $2^{r - C_6\log r + 1}$ strings of length at most $r - C_6\log r$, if $C_6$ is much larger than the implicit constant for $W$, then most elements of $W$ have to have complexity exceeding $r - C_6\log r$. Hence, essentially the same bound holds for $V$ as held for $W$:
  \begin{equation}\label{eq:boundOnV} 
  \vert V \vert \geq 2^{K^A_r(s) - O(\log r)}.
  \end{equation}

Now we can check that the conditions of Lemma \ref{lem:combinatorial} hold for $L$ and $V$. For all $v\in V\subseteq W$, by the definition of $W$ and \eqref{eq:LIsLarge},
\begin{align*}
    \vert N_v \vert &> 2^{K^A_r(\ell) - C_3 \log r}\\
    &\geq 2^{-(C_3 + 2C_1) \log r + 1} * 2^{K^A_r(\ell) + 2C_1 \log r + 1}\\
    &\geq 2^{-(C_3 + 2C_1) \log r + 1} * \vert L \vert. 
\end{align*}
Define $C_7= C_3 + 2 C_2 + 1$. Letting $\alpha = 2^{-C_7 \log r}$, we have
\begin{equation*}
    \vert N_v \vert \geq \alpha \vert L \vert,
\end{equation*}
for all $v\in V$. Hence, in order to apply Lemma \ref{lem:combinatorial}, it remains to define the relation $\sim$ between $u, v\in V$ and to show that most elements are not related to each other. 

Let $l$ be some line in $L$. If $u$ and $v$ are too close to each other, then $l(u)$ and $l(v)$ only determine $l$ to an unacceptably low precision. Additionally, we would like to guarantee that the pair $u, v$ has high complexity relative to $A$, as we will make use of this complexity in the final step of the proof. Hence, we say $u\sim v$ if either
\begin{enumerate}
    \item $\vert u - v\vert <2^{-t}$, or
    \item $K^A_r(u, v)\leq 2 r- C_8 \log r$ for some $C_8$ to be determined shortly.
\end{enumerate}

As in \cite{CholakCsorn2025Bourgain}, we know that for a given $v$, only $u$ that are in the same or an immediately adjacent $2^{-t}$ dyadic interval can be related, which implies
\begin{equation*}
    \{u \in V: \vert u - v\vert < 2^{-t}\}   \leq 3 * 2^{K^A_r(s) - K^A_t(s) - C_4 \log r}. 
\end{equation*}
Using \eqref{eq:boundOnV}, this becomes 
\begin{equation*}
    \{u \in V: \vert u - v\vert \leq 2^{-t}\}   \leq 2^{-K^A_t(s) - O(\log r)} \vert V\vert. 
\end{equation*}
Recalling that $s$ was random relative to $A$ and the definition of $t$, for $r$ sufficiently large,
\begin{align*}
    \{u \in V: \vert u - v\vert \leq 2^{-t}\}   &\leq 2^{-t - O(\log r)} \vert V\vert\\
    &= 2^{-\frac{\ve}{2n} r - O(\log r)} \vert V\vert\\
     &\leq 2^{-2C_7 \log r - 3} \vert V\vert\\
     &\leq \frac{\alpha^{2}}{8} \vert V\vert.
\end{align*}
Hence, by the definition of $\sim$, it remains only to show that 
\begin{equation}\label{eq:mostPointsComplex}
\{u \in V: K_r^A(u, v)\leq 2 r - C_8 \log r\}\leq \frac{\alpha^2}{8} \vert V \vert
\end{equation}
To see that typical pairs $(u, v)$ have high complexity relative to $A$, recall that $v$ satisfies $K_r^A(v)\geq r - C_6 \log r $ by our refinement of $W$. Hence,
\begin{align}
    K^A_r(u, v) &\geq K^A_r(v) + K^A_r(u \mid v) - O(\log r)\nonumber \\
    &\geq r + K^A_r(u \mid v) - O(\log r)\nonumber \\
    &\geq r + K^{A, v}_r(u) - O(\log r).\label{eq:oracleV}
\end{align}
Since
\begin{equation*}
    \vert \{u\in V: K^{A, v}_r(u) < r - \beta \log r\} \vert < 2^{r- \beta \log r +1},
\end{equation*}
applying \eqref{eq:boundOnV} and \eqref{eq:oracleV} ensures that there exists some $C_8$ such that \eqref{eq:mostPointsComplex} holds.

Having checked all its conditions, we apply Lemma \ref{lem:combinatorial} and conclude that there exist $u, v$ not related to each other such that
\begin{equation}\label{eq:largeIntersection}
    \vert N_u\cap N_v\vert \geq 2^{-2C_7\log r - 1}\vert L\vert.
\end{equation}
Fix this pair $u, v$ for the remainder of the argument.

\noindent \textbf{Obtaining the contradiction:} Combining \eqref{eq:LIsLarge}, which related the size of $L$ to the complexity of $\ell$, and \eqref{eq:largeIntersection}, we can pick some $l\in N_v \cap N_u$ such that
\begin{equation}\label{eq:surrogateComplexity}
    K_r^{A, u, v}(l) \geq K^A_r(\ell) - O(\log r). 
\end{equation}

We will complete the proof by bounding the quantity $K^A_r(l(u), l(v))$. On one hand, using the definition of $N_u$ and $N_v$,
\begin{align}
    K_r^A(l(u), l(v))&\leq K^A_r(l(u)) +  K_r^A(l(v)) + O(\log r)\nonumber \\
    &\leq 2\left( \frac{K^A_r(\ell)}{2} + r - \ve r\right)+ O(\log r)\nonumber \\
    &= K^A_r(\ell) + 2 r - 2 \ve r + O(\log r).\label{eq:finalUpperBound}
\end{align}

On the other hand, applying Lemma 22 of \cite{fiedler2025extensionsUnions} and the fact that $u$ and $v$ are separated by at least $2^{-t}$,
\begin{equation*}
    K_r^A(l(u), l(v)) \geq K^A_{r, r, r-t}(l(u), l(v), l).
\end{equation*}
We continue by reading off precision $r$ approximations of $u$ and $v$ from $l(u)$ and $l(v)$, applying the symmetry of information, and using the fact that the pair $(u, v)$ has high complexity,
\begin{align*}
    K^A_r(l(u), l(v))  &\geq  K^A_{r, r, r-t}(u, v, l) - O(\log r)\\
    &=  K^A_r(u, v) + K^A_{r-t, r, r}(l\mid u, v)- O(\log r)\\
    &\geq  K^A_r(u, v) + K^{A, u, v}_{r-t}(l)- O(\log r)\\
    &\geq 2 r + K^{A, u, v}_{r-t}(l)- O(\log r).\\ 
\end{align*}
Recall that \eqref{eq:surrogateComplexity} gave us a bound on the complexity of $l$ at precision $r$, not precision $r-t$. Our last step will be to use \eqref{eq:pseudoCaseLutz} and the definition of $t$ to bridge this gap.  
\begin{align}
    K^A_r(l(u), l(v)) &\geq 2 r + K^{A, u, v}_{r}(l) - K_{r, r-t}^{A, u, v}(l)  - O(\log r)\nonumber \\ 
    &\geq 2 r + K^{A, u, v}_{r}(l) - 2(n-1) t  - O(\log r)\nonumber \\
    &= 2 r + K^{A, u, v}_{r}(l) - 2(n-1) \frac{\ve}{2 n} r  - O(\log r)\nonumber \\
    &> 2 r + K^{A, u, v}_{r}(l) - \ve r  - O(\log r)\nonumber \\ 
    &\geq 2 r + K^A_r(\ell) - \ve r  - O(\log r).\label{eq:finalLowerBound} 
\end{align}
Combining \eqref{eq:finalUpperBound} and \eqref{eq:finalLowerBound} yields
\begin{equation*}
    \ve r \leq O(\log r),
\end{equation*}
 and completes the proof, since this is false for sufficiently large $r$.

\end{proofof}

\section{Acknowledgments}

I am grateful to Don Stull for his helpful comments on an earlier draft of this paper.

\bibliographystyle{amsplain}
\bibliography{references}

\providecommand{\bysame}{\leavevmode\hbox to3em{\hrulefill}\thinspace}
\providecommand{\MR}{\relax\ifhmode\unskip\space\fi MR }
\providecommand{\MRhref}[2]{%
  \href{http://www.ams.org/mathscinet-getitem?mr=#1}{#2}
}
\providecommand{\href}[2]{#2}
\begin{thebibliography}{10}

\bibitem{Bourgain10}
Jean Bourgain, \emph{The discretized sum-product and projection theorems}, J. Anal. Math. \textbf{112} (2010), 193--236. \MR{2763000}

\bibitem{bushling2025extension}
Ryan E.~G. Bushling and Jacob~B. Fiedler, \emph{Bounds on the dimension of lineal extensions}, J.~Fractal Geom. \textbf{12} (2025), no.~1/2, 105–--133.

\bibitem{case2015dimension}
Adam Case and Jack~H. Lutz, \emph{Mutual dimension}, ACM Trans.~Comput.~Theory \textbf{7} (2015), no.~3, 1--26.

\bibitem{CholakCsorn2025Bourgain}
Peter Cholak, Marianna Cs\"ornyei, Neil Lutz, Patrick Lutz, Elvira Mayordomo, and D.~M. Stull, \emph{Algorithmic information bounds for distances and orthogonal projections}, 2025.

\bibitem{csornyei2025improvedboundsradialprojections}
Marianna Cs\"ornyei and D.~M. Stull, \emph{Improved bounds for radial projections in the plane}, 2025.

\bibitem{downey2010}
Rodney~G. Downey and Denis~R. Hirschfeldt, \emph{Algorithmic randomness and complexity}, Springer, New York, 2010.

\bibitem{falconer2016strong}
Kenneth~J. Falconer and Pertti Mattila, \emph{Strong {M}arstrand theorems and dimensions of sets formed by subsets of hyperplanes}, J.~Fractal Geom. \textbf{3} (2016), no.~4, 319–--329.

\bibitem{fiedler2025extensionsUnions}
Jacob~B. Fiedler, \emph{On the packing dimension of unions and extensions of $k$-planes}, 2025.

\bibitem{fiedler2025universalsetsprojections}
Jacob~B. Fiedler and D.~M. Stull, \emph{Universal sets for projections}, 2025.

\bibitem{heraKeleti2019hausdorff}
Korn\'{e}lia H\'{e}ra, Tam\'{a}s Keleti, and Andr\'{a}s M\'{a}th\'{e}, \emph{Hausdorff dimension of unions of affine subspaces and of {F}urstenberg-type sets}, J.~Fractal Geom. \textbf{6} (2019), no.~3, 263--284.

\bibitem{Hera2019}
Kornélia Héra, \emph{Hausdorff dimension of {F}urstenberg-type sets associated to families of affine subspaces}, Annales Fennici Mathematici \textbf{44} (2019), no.~2, 903–923.

\bibitem{keleti2016lines}
Tam\'{a}s Keleti, \emph{Are lines much bigger than line segments?}, Proc.~Amer.~Math.~Soc. \textbf{144} (2016), no.~4, 1535--1541.

\bibitem{keleti2022equivalences}
Tamás Keleti and András Máthé, \emph{Equivalences between different forms of the {K}akeya conjecture and duality of {H}ausdorff and packing dimensions for additive complements}, 2023.

\bibitem{Lutz03a}
Jack~H. Lutz, \emph{Dimension in complexity classes}, {SIAM} J. Comput. \textbf{32} (2003), no.~5, 1236--1259.

\bibitem{Lutz03b}
\bysame, \emph{The dimensions of individual strings and sequences}, Inf. Comput. \textbf{187} (2003), no.~1, 49--79.

\bibitem{lutz2018algorithmic}
Jack~H. Lutz and Neil Lutz, \emph{Algorithmic information, plane {K}akeya sets, and conditional dimension}, ACM Trans.~Comput.~Theory \textbf{10} (2018), no.~2, 1--22.

\bibitem{LuLuMay2023PtS}
Jack~H. Lutz, Neil Lutz, and Elvira Mayordomo, \emph{Extending the reach of the point-to-set principle}, Information and Computation \textbf{294} (2023), 105078.

\bibitem{LutStu18Projections}
Neil Lutz and D.~M. Stull, \emph{Projection theorems using effective dimension}, 43rd International Symposium on Mathematical Foundations of Computer Science (MFCS 2018), 2018.

\bibitem{lutz2020bounding}
Neil Lutz and Donald~M. Stull, \emph{Bounding the dimension of points on a line}, Inform.~and Comput. \textbf{275} (2020), 104601.

\bibitem{Mayordomo02}
Elvira Mayordomo, \emph{A {K}olmogorov complexity characterization of constructive {H}ausdorff dimension}, Inf. Process. Lett. \textbf{84} (2002), no.~1, 1--3.

\bibitem{stull2022pinneddistancesetsusing}
D.~M. Stull, \emph{Pinned distance sets using effective dimension}, 2022.

\bibitem{wang2025volumeestimatesunionsconvex}
Hong Wang and Joshua Zahl, \emph{Volume estimates for unions of convex sets, and the {K}akeya set conjecture in three dimensions}, 2025.

\end{thebibliography}

 \end{document}